\documentclass{elsartArXiv}
\usepackage{array,amsmath,amstext,amssymb,amsbsy,graphicx,array,a4wide}
\usepackage{pdfsync}
\usepackage{bm}
\usepackage{subfig}
 \usepackage{color}
\usepackage{mathtools}
\usepackage{graphicx}
\usepackage{todonotes}

\renewenvironment{cases}{\begin{dcases}}{\end{dcases}}

\newcommand{\rd}[1]{\text{d}#1}

\newcommand{\bfv}{{\boldsymbol v}}
\newcommand{\bfC}{{\boldsymbol C}}
\newcommand{\bfA}{{\boldsymbol A}}

\newcommand{\bfsig}{{\boldsymbol\sigma}}

\newcommand{\bfmL}{\mathbf L}

\newcommand{\bfE}{{\boldsymbol E}}
\newcommand{\bfS}{{\mathbf S}}

\def\bfF{{\boldsymbol F}}
\def\bfPhi{{\boldsymbol\Phi}}

\newcommand{\bfX}{\boldsymbol X}

\newcommand{\bfN}{\boldsymbol N}
\newcommand{\bfJ}{\mathbf J}
\newcommand{\bfb}{{\mathbf b}}
\newcommand{\bfxi}{{\boldsymbol\xi}}
\newcommand{\bfa}{{\mathbf a}}
\newcommand{\bfB}{{\mathbf B}}
\newcommand{\bfD}{{\mathbf D}}

\newcommand{\bfK}{{\mathbf K}}

\newcommand{\bfI}{{\mathbf I}}

\newcommand{\bfu}{{\boldsymbol u}}

\newcommand{\bfw}{{\boldsymbol w}}
\newcommand{\bfn}{{\boldsymbol n}}

\newcommand{\bfeps}{{\boldsymbol\varepsilon}}
\newcommand{\bff}{{\boldsymbol f}}

\newcommand{\bfg}{{\boldsymbol g}}
\newcommand{\bfx}{{\boldsymbol x}}

\newcommand{\deriv}[2]{\frac{\partial #1}{\partial #2}}

\newcommand{\jac}[1]{\frac{\partial #1}{\partial\bfX}}

\newcommand{\bbR}{{\mathbb{R}}}

\begin{document}
\begin{frontmatter}

\title{Modified midpoint integration rule for the trilinear element in large deformation elasticity}
\author{Mirza Cenanovic,}  
\author{Peter Hansbo,}  
\author{David Samvin}  
\address{Department of Mechanical Engineering, J\"onk\"oping University,
SE-55111 J\"onk\"oping, Sweden}
\date{Received: date / Accepted: date}

\maketitle

\begin{abstract}
In this paper we
suggest two modified one-point  Gauss integration rules for the $Q_1$ bi-- or trilinear element.
The modifications both stabilize the hourglass modes of the one-point rule, and one of them is accurate also on severely distorted elements.
We
investigate the performance of the integration rules for the hexahedron element, and combine standard one--point integration of the volumetric terms with the modified 
rules for the isochoric terms to handle near incompressible situations.

\end{abstract}
\end{frontmatter}
\maketitle

\section{Introduction}

We consider the construction of an hourglass control method in the form of a modified Gaussian rule including derivative information 
to allow for one point integration of the full element matrices without loss of stability. 
The one point rule suffers from loss of stability leading to so-called hourglass modes. In order to remedy this, 
hourglass control methods have been proposed. In the context of large deformations, these are typically based on 
introducing auxiliary strain fields, cf. \cite{GlAr97,ReKuRe99,ReWr00}, or to use
Taylor expansions of the strain field, for linear problems by Schultz \cite{Sc85} and by Liu, Ong, and Uras \cite{Liu85}, with extension to large deformation problems in Liu, Belytschko, Ong, and Law \cite{LiBeOnLa85} and Reese, Wriggers, and Reddy \cite{ReWrRe00}. 
In contrast, we employ a modification of the one point Gauss rule using derivative information with selective differentiation,
an idea introduced and analyzed by Hansbo \cite{Ha98} for a linear Poisson problem, and subsequently used in combination with a non--confor\-ming affine tetrahedron element for a large deformation model in Hansbo and Larsson \cite{HaLa16}. For non-affine elements, there
is however a strong stiffening effect on slender structures when the element geometry is distorted, as shown in the numerical examples below (cf. also Reese et al. \cite{ReWrRe00}). 
The purpose of this paper is to modify the scheme in \cite{HaLa16} so that the sensitivity to element geometry is remedied.

An outline of the paper is as follows. In Section \ref{one} we consider the linear elasticity problem to introduce the basic ideas of our integration scheme; in Section \ref{two} we introduce our large deformation problem, which is formulated in terms of the first Piola--Kirchhoff stress corresponding to an isotropic Mooney--Rivlin model, and show how the numerical integration scheme is
implemented in this case. In Section \ref{numex} we present numerical examples and in Section \ref{four} we give some concluding remarks.  Implementation details are provided in an Appendix.

\section{The small deformation elasticity problem\label{one}}

\subsection{Continuous model}

We consider a elasticity problem on a domain $\Omega\subset\bbR^n$. For notational simplicity we shall set $n=3$, but 
the two--dimensional case will also be considered in the numerical examples.
We first consider the following linearized elasticity problem: 
Find the displacement $\bfu$ such that
\begin{align}
\bfsig    =  \lambda ~\nabla\cdot\bfu~\bfI 
   + 2 \mu \bfeps(\bfu) & \quad\text{in}\; \Omega, \\
   -\nabla \cdot \bfsig   = \bff &\quad\mbox{ in } \Omega \label{eq:pde}\\
\bfu ={\mathbf 0} &\quad\mbox{ on } \partial \Omega_\text{D}\\
\bfsig\cdot\bfn ={\bfg} &\quad\mbox{ on } \partial \Omega_\text{N}\\
\end{align}
Here $\lambda$ and $\mu$ are the 
Lam\'{e} parameters, which we assume constant in
$\Omega$. 
In terms of the modulus of elasticity, $\rm{E}$, and Poisson's ratio, $\nu$, we have
\[
\lambda = \frac{\rm{E}\,\nu}{(1-2\nu)(1+\nu)},\quad\mu=\frac{\rm{E}}{2\, (1+\nu)}.
\]
Furthermore,
$\bfeps\left(\bfu\right) = \left[\varepsilon^{ij}(\bfu)\right]_{i,j=1}^3$ 
is the strain tensor with components
\[ \varepsilon_{ij}(\bfu ) = \frac{1}{2}\left( \frac{\partial
  u_i}{\partial x_j}+\frac{\partial u_j}{\partial x_i}\right) ,
\] 
$\nabla\cdot\bfsig = \left[\sum_{j=1}^3\partial
  \sigma_{ij}/\partial x_j\right]_{i=1}^3$,
$\bfI = \left[\delta_{ij}\right]_{i,j=1}^3$ with $\delta_{ij} =1$
if $i=j$ and $\delta_{ij}= 0$ if $i\neq j$, and $\bfn$ is the
outward  normal to
$\Omega$ on the boundary $\partial\Omega = \partial\Omega_\text{D}\cup \partial\Omega_\text{N}$

On weak form the problem is to find $\bfu\in V$, 
\[
V=\{\bfv\in [H^1(\Omega)]^3:\; \bfv = {\mathbf 0}\; \text{on}\, \partial\Omega_\text{D}\} ,
\]
such that
\begin{equation}\label{eq:standard-weak}
a(\bfu,\bfv) = l(\bfv) 
\quad \forall \bfv \in V, 
\end{equation}
where 
\begin{equation}
a(\bfv,\bfw) =  2\mu \int_{\Omega}\bfeps(\bfv):\bfeps(\bfw)\,\text{d}\Omega 
+ \lambda \int_{\Omega}\nabla\cdot\bfv\,\nabla\cdot\bfw\,\text{d}\Omega 
\end{equation}
and 
\begin{align}
l(\bfv) &= \int_{\Omega}\bff\cdot\bfv\,\text{d}\Omega + \int_{\partial \Omega_N} \bfg\cdot\bfv\, \text{d}s .
\end{align}
\
\subsection{Finite element model}

Let $\mathcal{T}_h$ be a conforming mesh on $\Omega\subset\bbR^n$. On a reference brick $\hat{T}$ we have 8 basis functions $\hat{\varphi}_i(\bfxi)$, 
and 
we
define a map $\bfF_{T} : (\xi,\eta,\zeta)\rightarrow (x,y,z)$ from the reference configuration to a
physical element $T$ by
\[
\bfx=\bfF_{T}(\bfxi) := \sum_{i=1}^8\hat{\varphi}_i(\bfxi) \bfx_i ,
\]
where $\bfx_i$ is the physical location of element node $i$.
Under this map, we define $\varphi_i(\bfx) = \hat{\varphi}_i(\bfxi)$. Denoting $Q_1$ as the set of trilinear functions, our discrete space can then be written
\begin{equation}\label{spacev}\begin{array}{lll}
V_h &:= &\{ \bfv\in V: {\bfv\vert_T \in [Q_1(\bfF_{T}^{-1}T)]^3,\; \forall T\in\mathcal{T}_h}\},
\end{array}\end{equation}
where $\bfF_T^{-1}T$ denotes the map back to the reference element $\hat{T}$, and we can introduce the finite element method as: find $\bfu^h\in V^h$
such that
\begin{equation}\label{membraneFEM}
a(\bfu^h,\bfv) = l(\bfv) ,\quad \forall \bfv \in V^h .
\end{equation}
In the practical implementation we let
\[
\bfu^h = \left[\begin{array}{ccccccc}
\varphi_1 & 0 & 0 & \varphi_2 & 0 & 0  & \ldots \\[4mm]
0 & \varphi_1 & 0 & 0 & \varphi_2 & 0 & \ldots \\[4mm]
0 & 0 & \varphi_1 & 0 & 0  & \varphi_2 & \ldots \\[4mm]
\end{array}\right] \left[\begin{array}{c}
U_{x_1,1}\\[4mm]
U_{x_2,1}\\[4mm]
U_{x_3,1}\\[4mm]
U_{x_1,2}\\[4mm]
U_{x_2,2}\\[4mm]
U_{x_3,2}\\[4mm]
\vdots\end{array}\right] = \bfPhi\bfa
\]
where $(U_{x_1,i},U_{x_2,i},U_{x_3,i})$ is the approximate displacement in node number $i$.
We then arrive at the discrete matrix problem $\bfK\bfa = \bfb$
where $\bfb$ contains the load contributions and
\begin{equation}\label{stiffmat}
\bfK = \int_{\Omega} \bfB^{\rm T}\bfD\bfB \,\text{d}x_1\text{d}x_2\text{d}x_3
\end{equation}
where, with $\tilde\nabla$ the matrix differential operator
\[
\tilde\nabla := \left[\begin{array}{>{\displaystyle}c>{\displaystyle}c>{\displaystyle}c>{\displaystyle}c>{\displaystyle}c>{\displaystyle}c}
\frac{\partial}{\partial x_1} & 0 & 0 &\frac{\partial}{\partial x_2} & 0 & \frac{\partial}{\partial x_3} \\[4mm]
0 & \frac{\partial}{\partial x_2} & 0 & \frac{\partial}{\partial x_1} & \frac{\partial}{\partial x_3} & 0\\[4mm]
0 & 0 & \frac{\partial}{\partial x_3} & 0 &  \frac{\partial}{\partial x_1} & \frac{\partial}{\partial x_1}\end{array}\right]
\]
$\bfB = \tilde\nabla^{\rm T}\bfPhi$ and 
$\bfD$ is a matrix of the elastic moduli (cf. \cite{crisfield2012nonlinear}). The entries of $\bfB$ are typically computed using the inverse of the Jacobian mapping from a reference element,
which we shall assume is given on the domain $0\leq \xi_i\leq 1$.
We define
\[
\bfJ := \left[\begin{array}{>{\displaystyle}c>{\displaystyle}c>{\displaystyle}c}
\deriv{x_1}{\xi_1} & \deriv{x_2}{\xi_1} & \deriv{x_3}{\xi_1} \\[4mm]
\deriv{x_1}{\xi_2} & \deriv{x_2}{\xi_2} & \deriv{x_3}{\xi_2} \\[4mm]
\deriv{x_1}{\xi_3} & \deriv{x_2}{\xi_3} & \deriv{x_3}{\xi_3}
\end{array}\right] 
\] 
and then
\[
\left[\begin{array}{>{\displaystyle}c}
\deriv{\varphi_i}{x_1} \\[4mm]
\deriv{\varphi_i}{x_2} \\[4mm]
\deriv{\varphi_i}{x_3} 
\end{array}\right] =\bfJ^{-1}\left[\begin{array}{>{\displaystyle}c}
\deriv{\hat{\varphi}_i}{\xi_1} \\[4mm]
\deriv{\hat{\varphi}_i}{\xi_2}\\[4mm]
\deriv{\hat{\varphi}_i}{\xi_3} 
\end{array}\right] =: \bfJ^{-1}\nabla_\xi \hat{\varphi}
\]
and $\bfK$ is assembled from element contributions
\begin{equation}
\bfK_T :=  \int_{T} \bfB^{\rm T}\bfD\bfB\,  \,\text{d}x_1\text{d}x_2\text{d}x_3 = \int_{\hat{T}} \bfB^{\rm T}\bfD\bfB\, J \,\text{d}\xi_1\text{d}\xi_2\text{d}\xi_3 ,
\end{equation}
where $J:= \text{det}\,\bfJ$.

\subsection{Numerical integration}\label{numint1}

The one point Gaussian (midpoint) rule for the element-wise evaluation of (\ref{stiffmat}) is well known to result in a singular matrix, with singularity related to the {\em hourglass mode.}\/ 
This is unfortunate since the midpoint rule would be sufficiently accurate in terms of convergence of the method, and could be used were it not for the stability problem. 
For this reason, the use of $2\times 2\times 2$ Gauss points is usually recommended. However, one may alternatively use the following integration rule\begin{align*}
\int_{a}^b f(\xi)\, \text{d}\xi \approx{}& \int_{a}^b \left( f(\xi_m)+f'(\xi_m)(\xi-\xi_m)+\frac12 f''(\xi_m)(\xi-\xi_m)^2\right)\text{d}\xi \\
= {}& h\, f(\xi_m) + \frac{h^3}{24}f^{\prime\prime}(\xi_m) ,
\end{align*}
which is exact for polynomials up to second degree.
On our reference element we can use
the Taylor series expansion of  $f(\bfxi)$ around $\bfxi_m$,
\begin{align}\nonumber
f(\bfxi)\approx {} & f(\bfxi_m) + (\bfxi-\bfxi_m)\cdot\nabla_{\xi} f(\bfxi_m) \\
& + \frac12(\bfxi-\bfxi_m)\cdot\left((\bfxi-\bfxi_m)\cdot\nabla_\xi(\nabla_\xi f(\bfxi_m))\right).
\end{align}
Now, the second term in the expansion integrates to zero on the reference element by definition of the center of gravity, as do the mixed derivative terms in the third term.
We are left with
\begin{align}\nonumber
\int_{0}^1\int_{0}^1\int_{0}^1  f(\bfxi) \, \rd{\xi_1}\rd{\xi_2}\rd{\xi_3} \approx f(\bfxi_m) +\frac{1}{24}\sum_{i=1}^3\frac{\partial^2f}{\partial\xi_i^2}(\bfxi_m) .
\end{align}

We now want to apply this rule to the element--wise integration of contributions to $\bfK$.
We have
\begin{align}\nonumber
\frac{\partial^2}{\partial\xi^2_i}\left(\bfB^{\rm T}\bfD\bfB\, J\right) = {}&  \left(2\left(\frac{\partial\bfB}{\partial\xi_i}\right)^{\rm T}\bfD \frac{\partial\bfB}{\partial\xi_i} +\left(\frac{\partial^2\bfB}{\partial\xi^2_i}\right)^{\rm T}\bfD \bfB +\bfB^{\rm T}\bfD \frac{\partial^2\bfB}{\partial\xi^2_i}\right) J \\[3mm]
& + 2\left(\left(\frac{\partial\bfB}{\partial\xi_i}\right)^{\rm T}\bfD \bfB +\bfB^{\rm T}\bfD \frac{\partial\bfB}{\partial\xi_i}\right)\deriv{J}{\xi_i} +\bfB^{\rm T}\bfD\bfB\frac{\partial^2 J}{\partial \xi_i^2} ,\label{secondder}
\end{align}
and we thus seek suitable simplifications of this expression, notably avoiding second derivatives of involved quantities.

Defining
\[
\bfA_i := \bfJ^{-1}\deriv{\bfJ}{\xi_i},
\]
we have that
\[
\deriv{\bfJ^{-1}}{\xi_i} = -\bfA_i\bfJ^{-1}, 
\]
 and it follows that
\begin{equation}\label{simple}
\deriv{}{\xi_i}\left[\begin{array}{>{\displaystyle}c}
\deriv{\varphi_i}{x_1} \\[4mm]
\deriv{\varphi_i}{x_2} \\[4mm]
\deriv{\varphi_i}{x_3} 
\end{array}\right] =\bfJ^{-1}\left[\begin{array}{>{\displaystyle}c}
\frac{\partial^2\hat{\varphi}_i}{\xi_1\xi_i} \\[4mm]
\frac{\partial^2\hat{\varphi}_i}{\xi_2\xi_i}\\[4mm]
\frac{\partial^2\hat{\varphi}_i}{\xi_3\xi_i}
\end{array}\right] -\bfA_i\left[\begin{array}{>{\displaystyle}c}
\deriv{\varphi_i}{x_1} \\[4mm]
\deriv{\varphi_i}{x_2} \\[4mm]
\deriv{\varphi_i}{x_3} 
\end{array}\right] .
\end{equation}
These quantities are then used to build derivatives of the $\bfB$ matrices with respect to reference coordinates. 

{\em In the following, it is understood that all terms in the approximate integration schemes are evaluated at the midpoint of the reference element.}

We note that in case $\bfA_i$ is nonzero, there is an explicit dependence on $i$ in the derivative of $\bfB$. 
In the following we shall therefore write $\partial\bfB(\bfA_i)/\partial \xi_i$ for clarity. 
The simplest possible integration schemes using approximations of (\ref{secondder}) are
\begin{enumerate}
\item[A.] Approximate $\bfJ$ by a constant matrix. This leads to $\bfA_i = {\bf 0}$ and
\[
\frac{\partial^2}{\partial\xi^2_i}\left(\bfB^{\rm T}\bfD\bfB\, J\right) \approx 2\left(\frac{\partial\bfB(\bf{0})}{\partial\xi_i}\right)^{\rm T}\bfD \frac{\partial\bfB(\bf{0})}{\partial\xi_i} J .
\]
\item[B.] Approximate $\partial\bfJ^{-1}/\partial\xi_i$ by a constant matrix and $J$ by a constant. This leads to 
\[
\frac{\partial^2}{\partial\xi^2_i}\left(\bfB^{\rm T}\bfD\bfB\, J\right) \approx 2\left(\frac{\partial\bfB(\bfA_i)}{\partial\xi_i}\right)^{\rm T}\bfD \frac{\partial\bfB(\bfA_i)}{\partial\xi_i} J .
\]
\end{enumerate}
Our numerical experience is that including derivatives of $J$ does not improve the accuracy, so we will limit the discussion to these
two approximations.

The element matrix computed on the unit reference element $\hat{T}$ can hence be approximated as follows
\begin{align}
\bfK_T =  {}&\int_{\hat{T}} \bfB^{\rm T}\bfD\bfB\, J \,\text{d}\xi_1\text{d}\xi_2\text{d}\xi_3 \nonumber \\
\approx  {}& J \bfB^{\rm T}\bfD\bfB+
\frac{J}{12}\sum_{i=1}^3\left(\frac{\partial\bfB(\bfA_i)}{\partial\xi_i}\right)^{\rm T}\bfD\frac{\partial\bfB(\bfA_i)}{\partial\xi_i} .\label{scheme}
\end{align}
The first term in the approximation is the one obtained by a standard one point rule,  
and the remaining are called hourglass stabilization matrices. This integration scheme, with $\bfA_i={\mathbf 0}$, was suggested in \cite{Ha98,HaLa16}, and is related to classical hourglass control methods, cf. \cite{Sc85,Liu85,LiBeOnLa85}. Unfortunately, this method leads to overly stiff responses in cases of large deviation from affine mappings.
The same problem exists in other approaches to hourglass control, cf. Reese, Wriggers, and Reddy \cite{ReWrRe00}.
As we shall see in Section \ref{numex}, this is remedied by simply changing $\partial\bfB(\bf{0})/\partial \xi_i$ to $\partial\bfB(\bfA_i)/\partial \xi_i$ in (\ref{scheme}), which comes at the cost of a small number of additional $3\times 3$ matrix multiplications.

\section{A large deformation elasticity problem}\label{two}

\subsection{Continuous model}

Let an elastic body in its undeformed configuration, in a Cartesian coordinate system $\bfX$, occupy a three--dimensional domain $\Omega_0\subset{\mathbb{R}}^3$, with outward pointing normal $\bfN$ to the boundary $\partial\Omega_0$.
The corresponding position in the deformed domain
is denoted by $\bfx$, and the displacement field $\bfu$ is then given
by $\bfu(\bfX) = \bfx(\bfX)-\bfX$. 
Let us define the Jacobian of $\bff$ as
\begin{equation}
\frac{\partial \bff}{\partial \bfX}: = \left[\begin{array}{>{\displaystyle}l>{\displaystyle}l>{\displaystyle}l}
\frac{\partial f_1}{\partial X_1} & \frac{\partial f_1}{\partial X_2}& \frac{\partial f_1}{\partial X_3}\\[3mm]
\frac{\partial f_2}{\partial X_1} & \frac{\partial f_2}{\partial X_2}& \frac{\partial f_2}{\partial X_3}\\[3mm]
\frac{\partial f_3}{\partial X_1} & \frac{\partial f_3}{\partial X_2}& \frac{\partial f_3}{\partial X_3}
\end{array}\right]
.
\end{equation}
The deformation gradient on $\Omega_0$ is then given by
\begin{equation}
\bfF(\bfu) := \jac{\bfx}=\bfI + \jac{\bfu} 
\end{equation}
Consider next a potential energy functional given by
\begin{equation}
\Pi(\bfu) := \Psi(\bfu)-\Pi^\text{ext}(\bfu)
\end{equation}
where $\Psi$ is the strain energy functional and $\Pi^\text{ext}$ is the potential of external loads. We will assume conservative loading so that $\Pi^\text{ext}(\bfu)  = l(\bfu)$ is linear.
We have
\begin{equation}
\Psi(\bfu) = \int_{\Omega_0}\hat{\Psi}_{\bfX}(\bfE(\bfu))\, d\Omega_0 = \int_{\Omega_0}{\Psi}_{\bfX}(\bfC(\bfu))\, d\Omega_0 ,
\end{equation}
where $\Psi_{\bfX}$ is the strain energy per unit volume.
Assuming for simplicity that $\bfu$ is zero on part of the boundary $\partial\Omega_0$, then minimizing the potential energy leads to the variational problem of finding $\bfu: \Omega \rightarrow \bbR^3$ such that
\begin{equation}
\int_{\Omega_0} \bfS(\bfu) : {\bfE'}(\bfu,\bfv)\, d\Omega = l(\bfv) ,
\end{equation}
or
\[
A(\bfu,\bfv) = l(\bfv) ,
\]
for all $\bfv: \Omega_0 \rightarrow \bbR^3$ vanishing on $\partial\Omega_0$, where
\begin{equation}
 \bfS(\bfu)  := \frac{\partial\hat{\Psi}_{\bfX}}{\partial\bfE}=2\deriv{\Psi_{\bfX}}{\bfC}
\end{equation}
is the second Piola--Kirchhoff stress tensor and
\begin{equation}
{\bfE'}(\bfu,\bfv) := \frac12 \left(\bfF(\bfu)^{\rm T}\jac{\bfv} + \left(\jac{\bfv}\right)^{\rm T}\bfF(\bfu)\right)
\end{equation}
is the variational derivative of $\bfE$.

For definiteness, we use an isotropic Moo\-ney--Rivlin model in which we choose
parameters $E$ and $\nu$, and define
$K=E/(3(1-2\nu))$, $\mu=E/(2(1+\nu))$, and $K_1=K_2 = \mu/2$. Then the Mooney--Rivlin strain energy density is given by 
\begin{equation}
\Psi_{\bfX}(\bfC) :=  K_1  (\hat{I}_1-3) + K_2 (\hat{I}_2 -3) + \frac12 K (J_C-1)^2
\end{equation}
where $J_C := \text{det}\, \bfC$, $\hat{I}_1 := J_C^{-1/3}I_1$, and $\hat{I}_2 := J_C^{-2/3}I_2$, with $I_1$ and $I_2$ the first and second invariants of $\bfC$, cf. \cite{crisfield2012nonlinear}.

We can now introduce the finite element method: Find $\bfu^h\in V^h$, where
\begin{align}\nonumber
 V^h = {}& \{\bfv: \bfv\in [W^h_k]^3, \; \text{$\bfv$ is zero on $\partial\Omega_0$}\},
\end{align}
such that
\begin{equation}\label{membraneFEM}
A(\bfu^h,\bfv) = l_{h}(\bfv) ,\quad \forall \bfv \in V^h,
\end{equation}
For the solution of the nonlinear problem (\ref{membraneFEM}) we adopt Newton iterations:
we compute updates $\Delta \bfu^h \in V^h$, for each previous iteration $\bfu^{h(k)}$, such that
\begin{equation}\label{membraneNewton}
A'(\bfu^{h(k)},\bfv,\Delta \bfu_h) = l_{h}(\bfv) - A(\bfu^{h(k)},\bfv),\; \forall \bfv \in V^h,
\end{equation}
resulting in the iterative solutions 
\[
\bfu^{h(k+1)}=\bfu^{h(k)}+\Delta \bfu^h \rightarrow \bfu^h\;\text{with $k$}. 
\]
In (\ref{membraneNewton}) we introduced the tangent form, or directional derivative, of $A(\cdot,\cdot)$, 
\begin{align}\nonumber
A'(\bfu,\bfv,\bfw) := &{} \int_{\Omega_0}  \bfE'(\bfu,\bfv) :\mathcal{L}:\bfE'(\bfu,\bfw)\, d\Omega \\
&{} + \int_{\Omega_0}\bfS(\bfu) :\left(\left(\jac{\bfv}\right)^{\rm T}\jac{\bfw}\right)\, d\Omega,\label{aprime}
\end{align}
where 
we used the material elasticity tensor defined as
\begin{equation}\label{consistentlin}
\mathcal{L} :=  \frac{\partial\bfS}{\partial\bfE}=2 \frac{\partial\bfS}{\partial\bfC}.
\end{equation}
This linearization can be split into a volumetric part and an isochoric part which will be used for the purpose of underintegrating the volumetric part to avoid locking in near incompressibility, cf. Appendix.

Following \cite{HaLa16} we propose to apply the integration rule to (\ref{membraneNewton}), after linearization.
This leads to a simple perturbed Newton method which in our computational experience shows similar convergence to Newton's method applied to a fully integrated element.
To define this approach, we employ Voigt notation following \cite{crisfield2012nonlinear} and
the element tangent stiffness matrix $\bfK_T$ can be computed as follows:
\begin{equation}
\bfK_T =\int_{T}\left( \bfB^{\text{T}}_L\bfmL\bfB_L +\bfB^{\text{T}}_{NL}{\mathbf T}\bfB_{NL}\right)\,\rd{x_1}\rd{x_2}\rd{x_3} \label{tangentStiffness}.
\end{equation}
Here $\mathbf T$ contains three copies of the second Piola--Kirchhoff tensor $\bfS$ on the diagonal, $\bfmL$ denotes the matrix representation of the fourth order tensor $\mathcal{L}$,
$\bfB_L$ is the Voigt representation of $\bfE'$ applied to the basis functions and $\bfB_{NL}$ is the Voigt representation of the Jacobian of the basis functions, see \cite{crisfield2012nonlinear} for details.
We now apply the one point integration formula to obtain
\begin{align}
\bfK_T \approx {}& J\left( \bfB^{\text{T}}_{L}\bfmL\bfB_{L} +\bfB^{\text{T}}_{NL}{\mathbf T}\bfB_{NL}\right) \nonumber\\
{}& +\frac{J}{12}\sum_{i=1}^3\left(\frac{\partial\bfB_{L}(\bfA_i)}{\partial\xi_i}\right)^{\rm T}\bfmL\frac{\partial\bfB_{L}(\bfA_i)}{\partial\xi_i} \nonumber\\
{}& +\frac{J}{12}\sum_{i=1}^3\left(\frac{\partial\bfB_{NL}(\bfA_i)}{\partial\xi_i}\right)^{\rm T}{\mathbf T}\frac{\partial\bfB_{NL}(\bfA_i)}{\partial\xi_i} \label{onep1}.
\end{align}

For the integration of $a_{h}(\bfu^{h(k)},\bfv)$ in the right hand side of the Newton iterations, 
we note that taking the $\xi_i$--derivative of $\bfS(\bfu) : {\bfE'}(\bfu,\bfv)$ gives rise to the same terms as the linearization 
and we we can approximate the contribution to the right-hand side as follows:
\begin{align}
\int_{T} \bfB^{\text{T}}_L\mathbf S  \,\rd{x_1}\rd{x_2}\rd{x_3}   \approx{}& J \bfB^{\text{T}}_{L}\mathbf S +\frac{ J}{12}\sum_{j=1}^3\frac{\partial\bfB^{\text{T}}_L(\bfA_i)}{\partial \xi_j}\bfmL\frac{\partial\bfB_L(\bfA_i)}{\partial\xi_j}{\mathbf a}^k \nonumber\\
{}& +\frac{ J}{12}\sum_{j=1}^3\frac{\partial\bfB^{\text{T}}_{NL}(\bfA_i)}{\partial \xi_j}{\mathbf T}\frac{\partial\bfB_{NL}(\bfA_i)}{\partial\xi_j}{\mathbf a}^k ,\label{onep2}
\end{align}
where ${\mathbf a}^k$ denotes the element displacement vector at iteration $k$ and $\mathbf S$  denotes the second Piola-Kirchhoff stress tensor on Voigt form.

As is evident from (\ref{onep1}) and (\ref{onep2}), this approximate Newton method only requires linearized quantities that are computed in standard Newton solvers for the fully integrated element.



\section{Numerical examples}\label{numex}
In the numerical examples below we use ``full integration'' to mean
an integration rule which is exact for polynomials up to and including
degree 2.

\subsection{Split into volumetric and isochoric parts}

In order to avoid locking in near incompressibility, we employ underintegration
of the volumetric part of $\mathbf{L}$ and $\mathbf{S}$ for the
case of Mooney-Rivlin material and define on Voigt form

\[
\begin{cases}
\mathbf{S}^{\text{iso}} & =\frac{\partial}{\partial{\bf C}}\left(K_{1}(\hat{I}_{1}-3)+K_{2}(\hat{I}_{2}-3)\right)\\
{\bf S}^{\text{vol}} & =\frac{\partial}{\partial{\bf C}}\left(\frac{1}{2}K(J_C-1)^{2}\right)\\
{\bf L}^{\text{iso}} & =\frac{\partial}{\partial{\bf C}}\mathbf{S}^{\text{iso}}\\
{\bf L}^{\text{vol}} & =\frac{\partial}{\partial{\bf C}}{\bf S}^{\text{vol}},
\end{cases}
\]
see Appendix for details on the implementation of the Voigt form.

\subsection{Underintegration}

We compare full integration of the volumetric and isochoric parts
of the discrete second Piola-Kirchhoff stress $\bm{T}$ and tangent
stiffness $\bm{L}$ to the full integration of the isochoric parts
$\bm{S}^{\text{iso}}$ and $\bm{L}^{\text{iso}}$ but one point integration
of the volumetric parts $\bm{T}^{\text{vol}}$ and $\bm{L}^{\text{vol}}$,
we denote this classical approach by ``1-PVol''. 

For our approach, we denote the full stabilization term as

\[
s_{h}:=\frac{J}{12}\sum_{j=1}^{3}\left(\frac{\partial\mathbf{B}_{L}^{\rm T}(\bfA_i)}{\partial\xi_{j}}({\bf L}^{\text{iso}}+{\bf L}^{\text{vol}})\frac{\partial\mathbf{B}_{L}(\bfA_i)}{\partial\xi_{j}}+\frac{\partial\mathbf{B}_{NL}^{\rm T}(\bfA_i)}{\partial\xi_{j}}({\bf T}^{\text{iso}}+{\bf T}^{\text{vol}})\frac{\partial\mathbf{B}_{NL}(\bfA_i)}{\partial\xi_{j}}\right)
\]
and the isochoric stabilization term as

\[
s_{h}^{\text{iso}}:=\frac{J}{12}\sum_{j=1}^{3}\left(\frac{\partial\mathbf{B}_{L}^{\rm T}(\bfA_i)}{\partial\xi_{j}}{\bf L}^{\text{iso}}\frac{\partial\mathbf{B}_{L}(\bfA_i)}{\partial\xi_{j}}+\frac{\partial\mathbf{B}_{NL}^{\rm T}(\bfA_i)}{\partial\xi_{j}}{\bf T}^{\text{iso}}\frac{\partial\mathbf{B}_{NL}(\bfA_i)}{\partial\xi_{j}}\right)
\]
 and thus we define ``1-PStab''
to mean 

\[
\mathbf{K}_{T}\approx J\left(\mathbf{B}_{L}^{\rm T}({\bf L}^{\text{iso}}+{\bf L}^{\text{vol}})\mathbf{B}_{L}+\mathbf{B}_{NL}^{\rm T}({\bf L}^{\text{iso}}+{\bf L}^{\text{vol}})\mathbf{B}_{NL}\right)+s_{h}
\]
and

\[
\mathbf{g}=\int_{T}\mathbf{B}_{L}^{\rm T}\mathbf{S}dx_{1}dx_{2}dx_{3}\approx J\mathbf{B}_{L}^{\rm T}({\bf S}^{\text{iso}}+{\bf S}^{\text{vol}})+s_{h}.
\]

Furthermore we define ``1-PStabIso'' to mean

\[
\mathbf{K}_{T}\approx J\left(\mathbf{B}_{L}^{\rm T}({\bf L}^{\text{iso}}+{\bf L}^{\text{vol}})\mathbf{B}_{L}+\mathbf{B}_{NL}^{\rm T}({\bf L}^{\text{iso}}+{\bf L}^{\text{vol}})\mathbf{B}_{NL}\right)+s_{h}^{\text{iso}}
\]
and

\[
\mathbf{g}=\int_{T}\mathbf{B}_{L}^{\rm T}\mathbf{S}dx_{1}dx_{2}dx_{3}\approx J\mathbf{B}_{L}^{\rm T}({\bf S}^{\text{iso}}+{\bf S}^{\text{vol}})+s_{h}^{\text{iso}}.
\]

Lastly, we consider the case when the stabilization term is defined using $\bfA_i=\mathbf{0}$, i.e., using a constant approximation for the Jacobian. We define``1-PStabConstJ" to mean $s_h|_{\bfA_i=\mathbf{0}}$ and ``1-PStabIsoConstJ" to mean $s_h^{\text{iso}}|_{\bfA_i=\mathbf{0}}$.

\subsection{Cantilever beam}

We consider the case of a cantilever with dimensions

\[
0\leq x_1\leq1/2,\ 0\leq x_2\leq1/10,\ 0\leq x_3\leq1/10,
\]

in meters and with zero displacements at $x=0$. This Cantilever was
subjected to an external virtual work

\[
l(\bm{v})=\int_{\Omega_0}\bm{f}\cdot\bm{v}d\Omega
\]

with $\bm{f}=(0,0,-10)\text{GN/m}^{3}$. The moduli of elasticity
were chosen as $E=200$GPa, $\nu=0.33$. For the Mooney-Rivlin model,
we chose to define the bulk moduls as $K=E/(3(1-\nu))$, $K_{1}=E/(2(1+\nu))$
and $K_{2}=0$ effectively reducing the model to a Neo-Hookean material
model. 

\begin{figure}
\begin{centering}
\subfloat[Solution using full integration. Maximum tip displacement is -0.3035.\label{fig:SolutionFullIntegration}]{\centering{}\includegraphics[width=0.4\columnwidth]{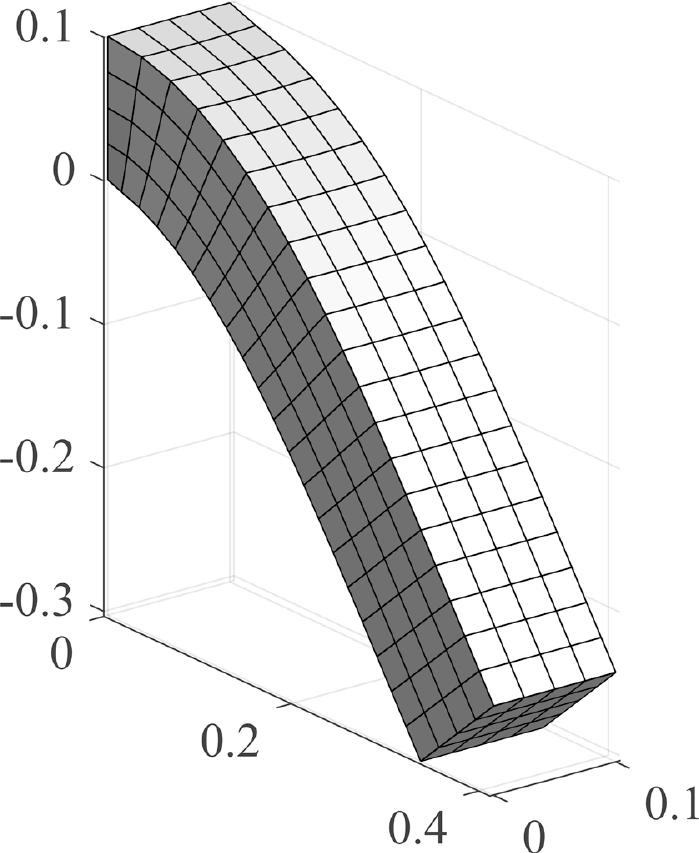}}\hspace{0.05\columnwidth}\subfloat[Solution using one point integration with stabilization. Maximum tip
displacement is -0.3037.\label{fig:SolutionOnePointIntegration}]{\begin{centering}
\includegraphics[width=0.4\columnwidth]{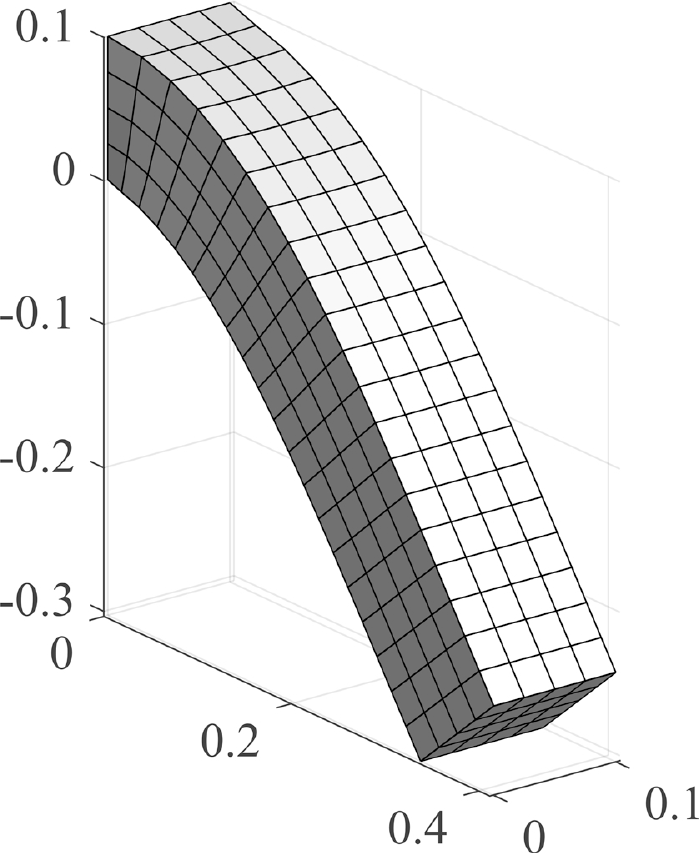}
\par\end{centering}
}\caption{Cantilever solution}
\par\end{centering}
\end{figure}

We compare the one point integration rule with stabilization, which
converges rapidly to the solution in Figure \ref{fig:SolutionOnePointIntegration},
to the result obtained using full integration in Figure \ref{fig:SolutionFullIntegration}
which shows that the difference in maximum tip displacement is small
for this particular meshsize.

\subsection{Stiffness comparison}

\begin{figure}
\centering{}\includegraphics[width=0.8\columnwidth]{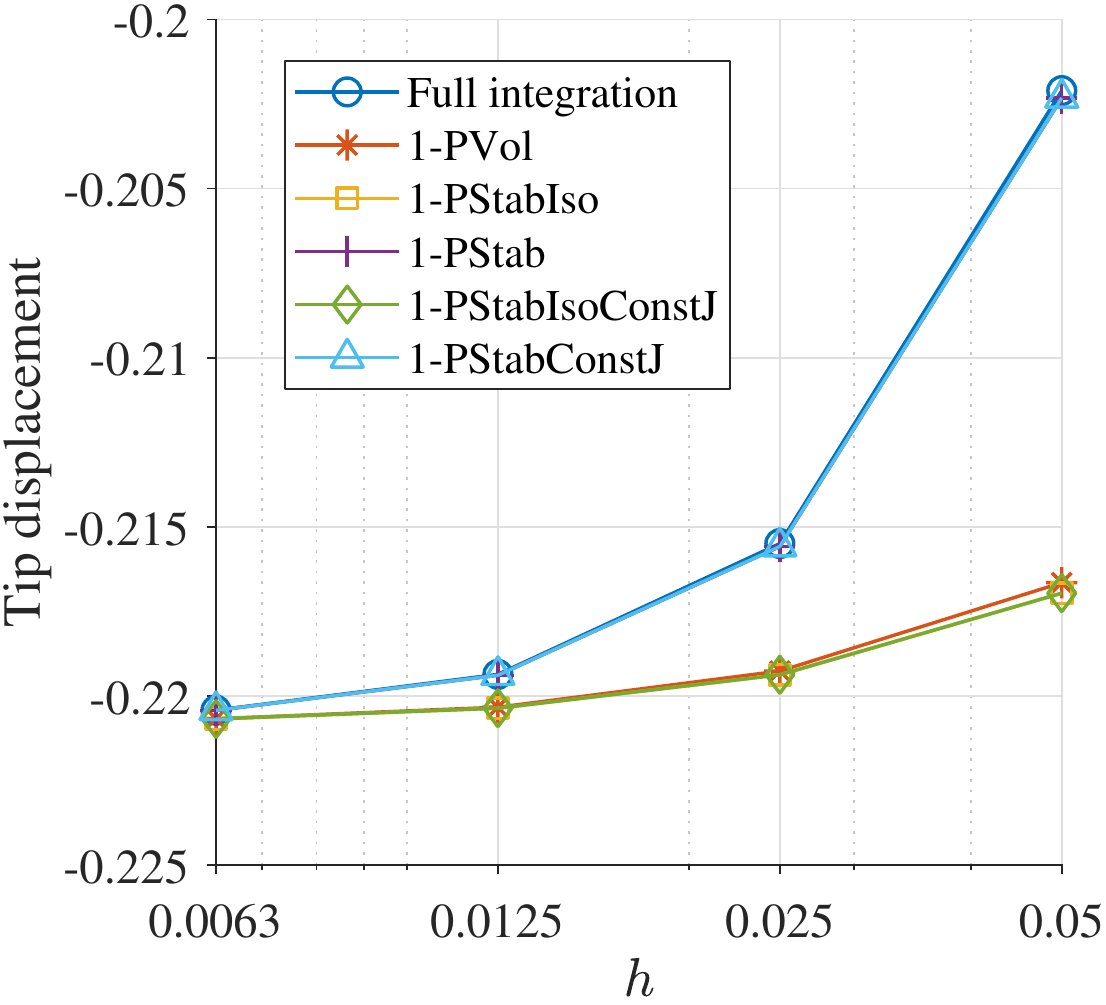}\caption{Tip displacement under fixed load and uniformly changing meshsize
using different integration methods.\label{fig:displacement-h}}
\end{figure}

We compare the tip deflection of the cantilever computed using different
integration schemes for a fixed load with decreasing meshsizes. We
compare the one-point integration with stabilization to full integration
and the classical reduced integration. Using the same parameters as
in the previous section we measure the tip displacement using different
meshsizes. We start with a uniform mesh with 2$\times$2$\times$10 elements and use
uniform refinement for each new mesh size. We show the obtained results
in Figure \ref{fig:displacement-h}, where it can be noticed that
the one point integration with stabilization behaves almost exactly
like the full integration and that the one point integration with
stabilized isochoric term gives nearly identical results as the classical
approach of one point integration of the volumetric term. The lack
of stabilization on the volumetric term makes the one point integration
with stabilization considerably softer than using stabilization on
both terms.

\begin{figure}
\centering{}\includegraphics[width=0.6\columnwidth]{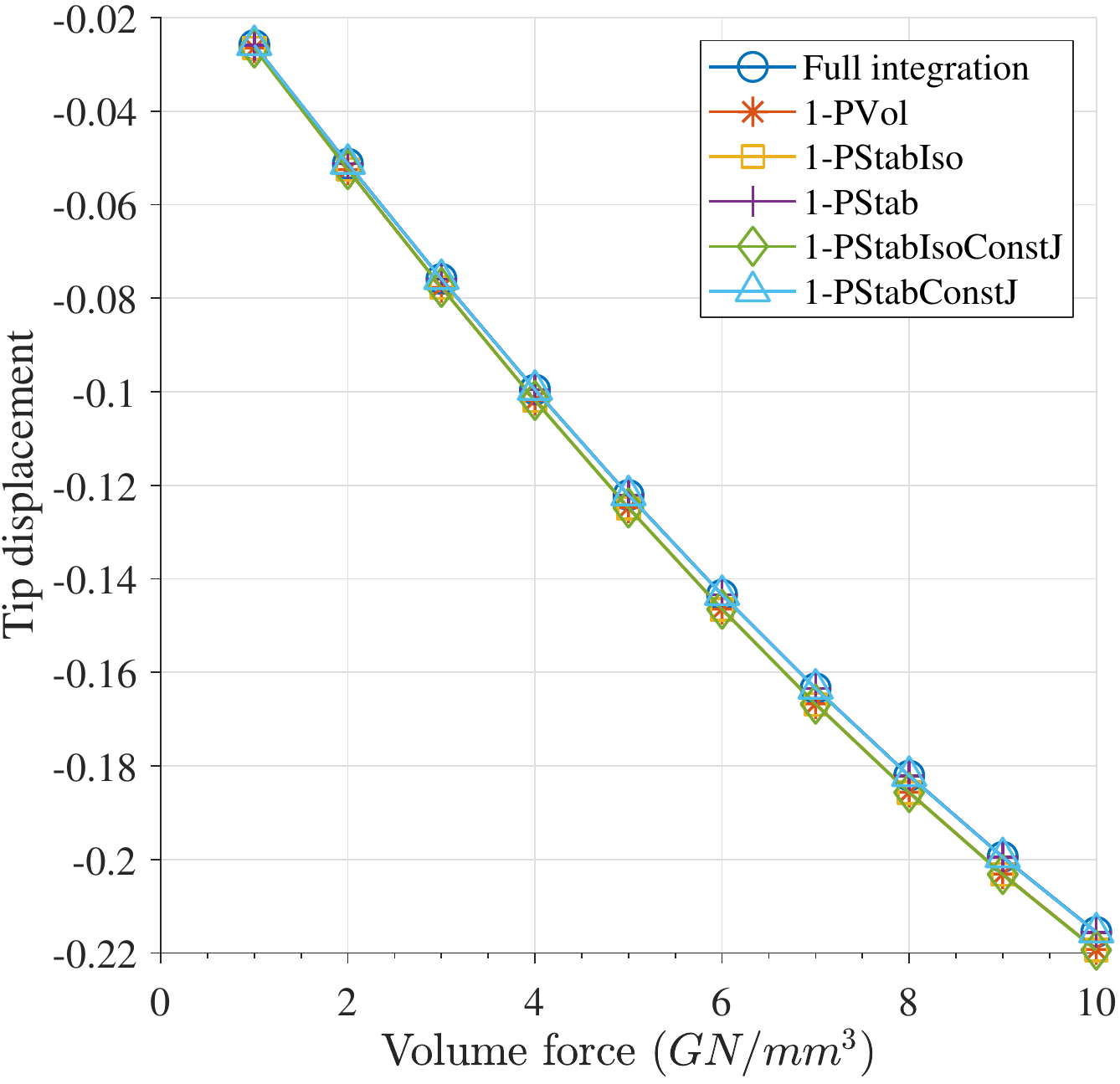}\caption{Tip displacement under varying load using different integration methods.\label{fig:displacement-f}}
\end{figure}

In Figure \ref{fig:displacement-f}, we compare the tip displacement
for the different integration approaches under linearly increasing
volume load. The load is linearly increasing from $\bm{f}=(0,0,-1)$
$\text{GN/m}^{3}$ to $\bm{f}=(0,0,-10)$ $\text{GN/m}^{3}$and it
can be seen that the one point integration with stabilization approach
continues to overlap with the full integration approach and that the
same overlapping can be seen for the one point integration with stabilization
of the isochoric term and the reduced integration.

\subsection{Near incompressibility}

\begin{figure}
\centering{}\includegraphics[width=0.7\columnwidth]{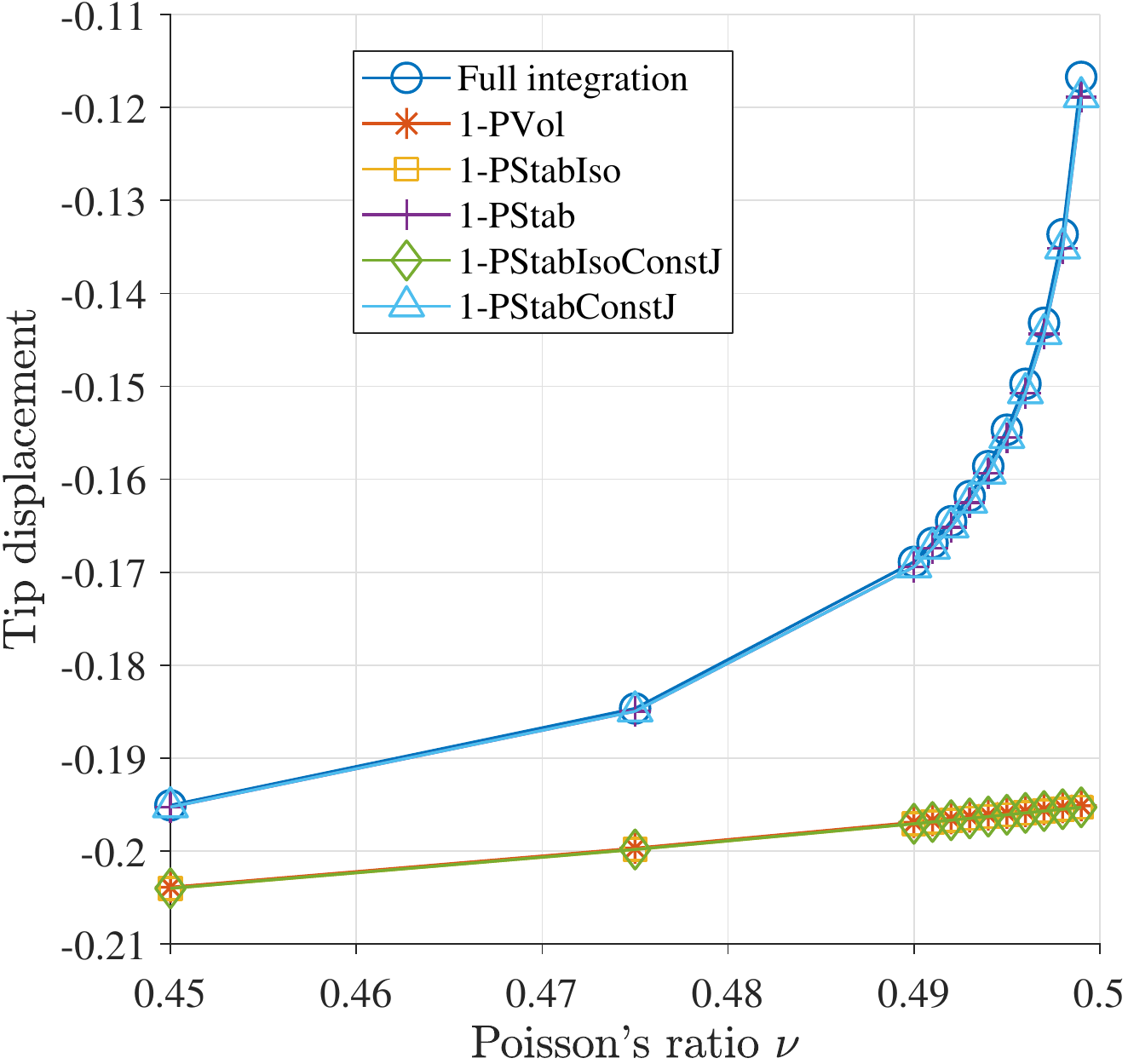}\caption{Locking behavior with Full integration / one point integration with
stabilization versus the locking free solution using reduced integration
/ one point integration of the isochoric term.\label{fig:displacement-nu}}
\end{figure}

In Figure \ref{fig:displacement-nu}, we show the tip displacement
of the same cantilever, load and material parameters from the previous
Section as $\nu\rightarrow1/2$. Notice the locking behavior when
using full integration as well as one point integration with stabilization
compared to one point integration of the volumetric term and one point
integration with stabilization of the isochoric term. In our numerical
test we were able to get stable solutions with $\nu$ reaching up
to $\nu=0.499$. 

\subsection{Mesh distortion\label{subsec:Mesh-distortion}}

\begin{figure}
\begin{centering}
\subfloat[Sketch of the beam.\label{fig:beamDistortSketch}]{\centering{}\includegraphics[width=0.4\columnwidth]{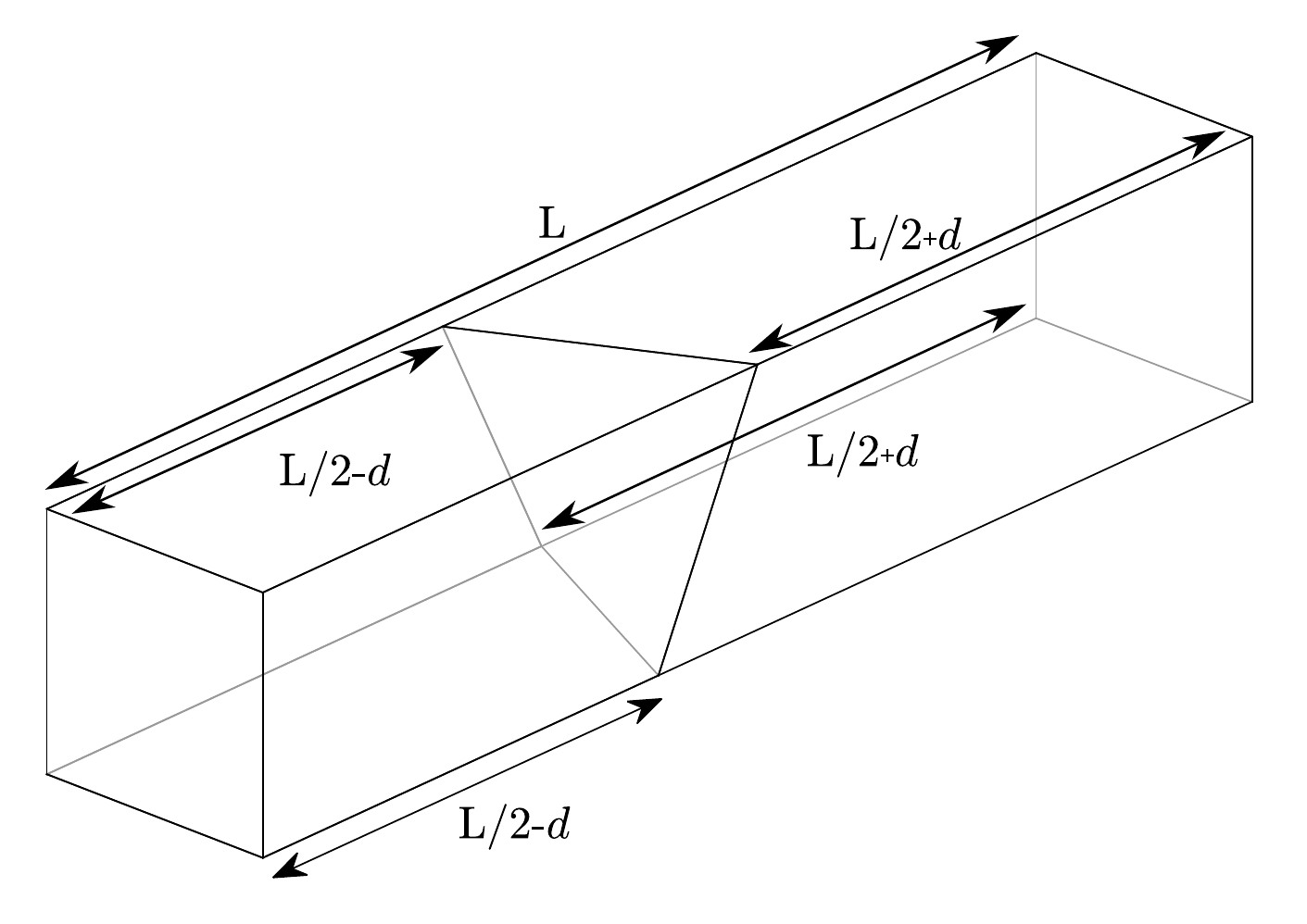}}\subfloat[Mesh of the beam.\label{fig:beamDistortMesh}]{\centering{}\includegraphics[width=0.5\columnwidth]{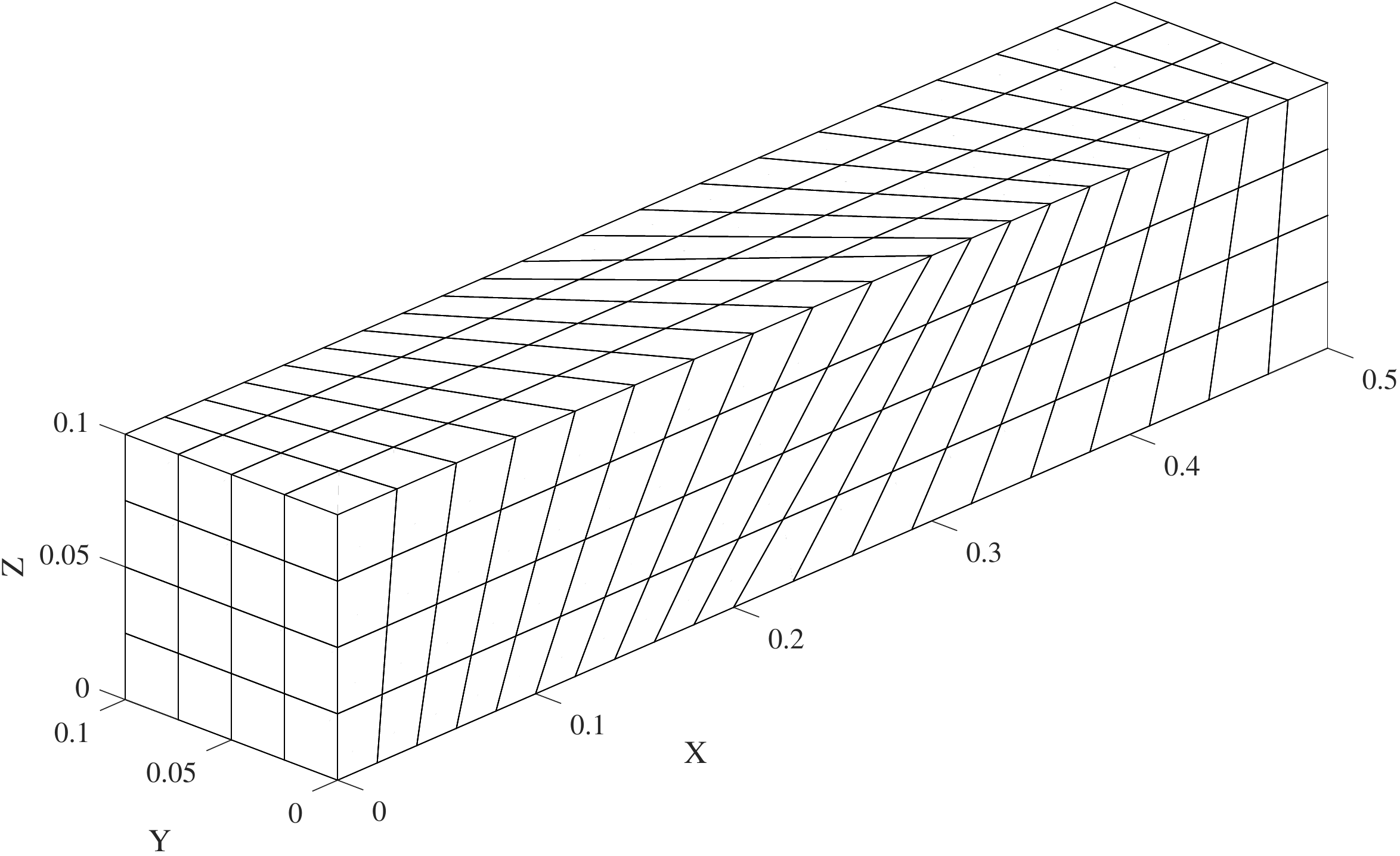}}\caption{Beam used to investigate sensitivity to mesh distortion}
\par\end{centering}
\end{figure}

In order to investigate the sensitivity of the proposed one point
integration schemes to mesh distortion we use a modified cantilever
beam with a distortion parameter $d$ as shown in Figure \ref{fig:beamDistortSketch}
and the corresponding mesh is shown i Figure \ref{fig:beamDistortMesh}.
All other parameters of the model are the same as in the previous
sections. We compared the relative change of the resultant of the
tip displacement $U_r$ with an increasing distortion parameter $d$ from
an initial distortion of $d=0$ to $d=0.2$, so that $U_r =U_r(d)$. The results are presented in
Figure \ref{fig:HB_Ur_d}, where the relative tip displacement is defined as
\[
U_{r}^{\text{rel}}=1-U_{r}(d)/U_{r}(0).
\]
Notice how the one point integration schemes
diverge from the reference tip displacement with almost 60 to above
70\% as the distortion becomes severe while the full integration and
the classical reduced integration scheme remain relatively unaffected. This sensitivity to mesh distorsion is compatible with the results of Reese et al. \cite{ReWrRe00}.

\begin{center}
\begin{figure}
\centering{}\includegraphics[width=0.6\columnwidth]{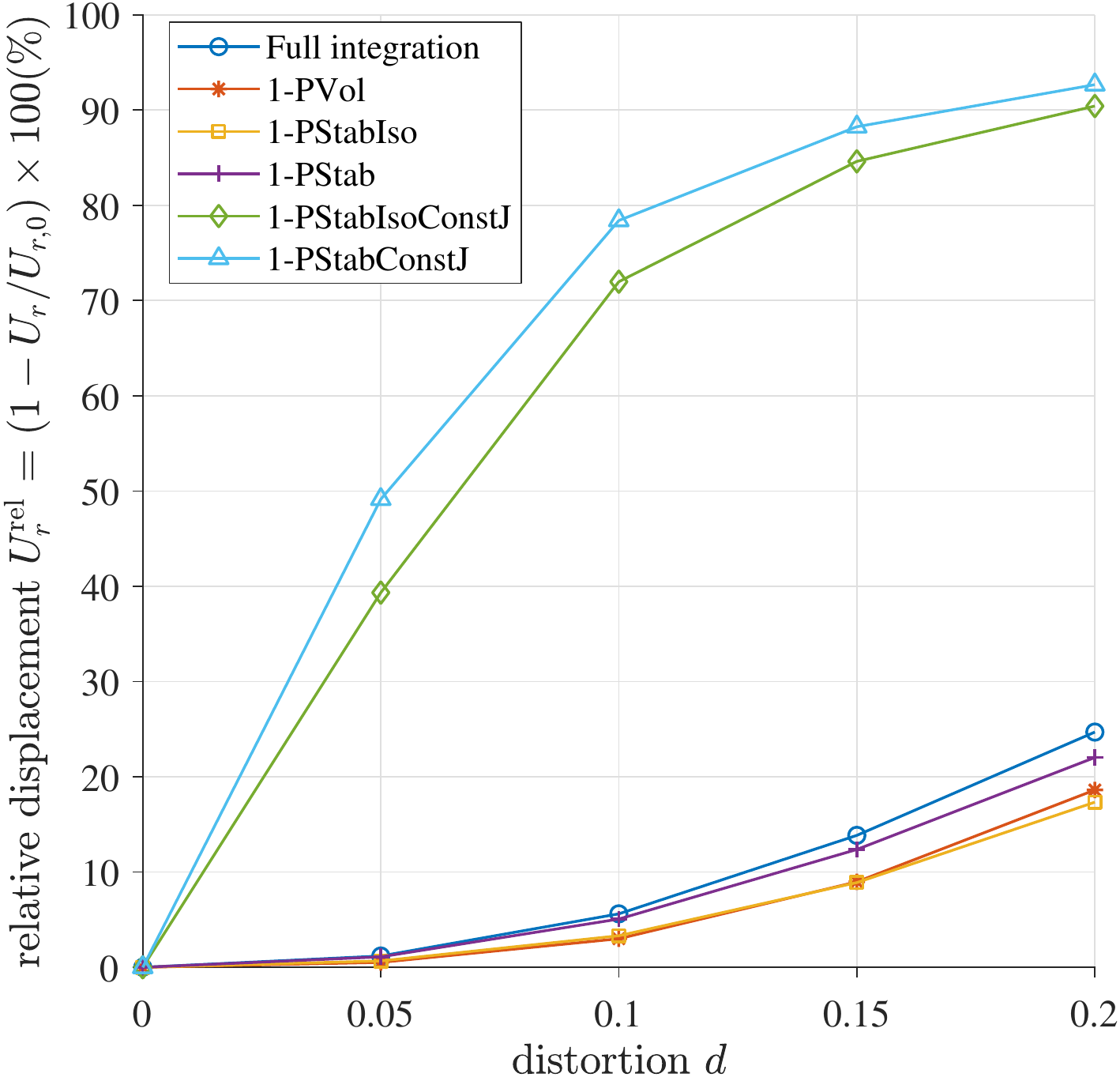}\caption{Effects of mesh distortion on the relative tip displacement for the non-linear problem.\label{fig:HB_Ur_d}}
\end{figure}
\par\end{center}

%
%
%

\subsection{Linear elastic mesh distortion}

We compare the sensitivity of the non-linear solution due to mesh
distortion from Section \ref{subsec:Mesh-distortion} to a linear
solution with the volumetric load $\bm{f}=(0,0,-1)$ $\text{GN/m}^{3}$.
The same beam geometry is used as in the previous section with the
same distortion and material parameters on the same mesh. The results
can be seen in Figure \ref{fig:HB_Ur_linear_d}, where a clear similarity
can be seen to the non-linear case, the one-point integration scheme
does not integrate the determinant of the Jacobian accurately enough.

\begin{center}
\begin{figure}
\centering{}\includegraphics[width=0.6\columnwidth]{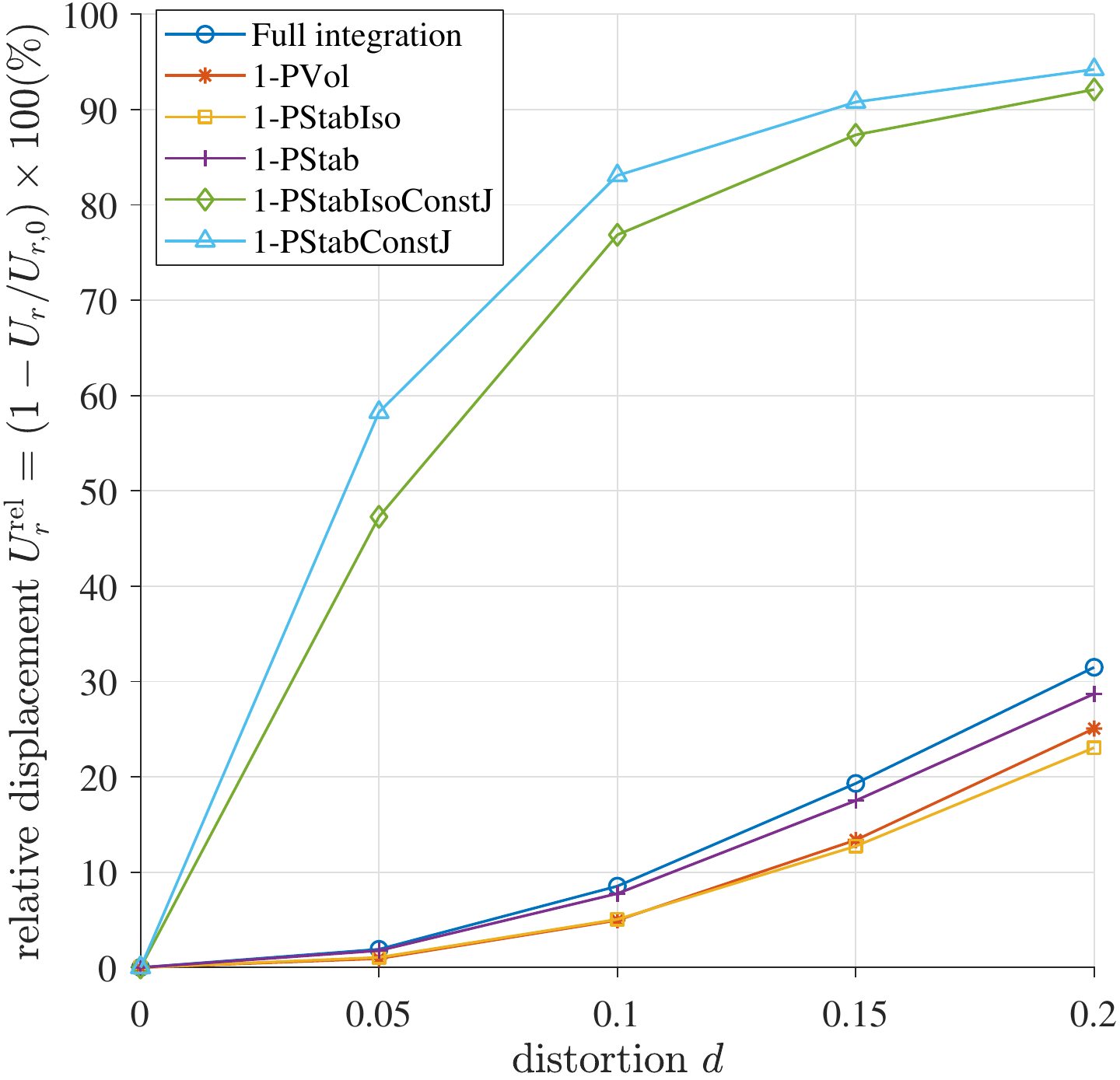}\caption{Effects of mesh distortion on the relative tip displacement for the linear problem.\label{fig:HB_Ur_linear_d}}
\end{figure}
\par\end{center}

\section{Concluding remarks\label{four}}

We have presented a general approach for one point integration of bi- and tri-linear finite elements in large deformation elasticity. The approach is
related to Taylor expansion ideas \cite{Sc85,ReWrRe00} but with a focus on integration.   

The basic method is highly accurate on meshes that do not contain highly distorted elements but shows an artificial stiffening as the element distortion increases, which is in accord with other hourglass stabilization methods. We therefore proposed a further modification which takes into account the effect of element distortion more accurately. This rule shows results close to the standard $2\times 2\times 2$ Gauss rule also on severely distorted meshes. 

\clearpage{}

\appendix
\section*{Appendix}\label{sec:Appendix}

\setlength{\jot}{10pt}

The linearisation is implemented using the Voigt form
and the details are given below. We follow the notation of de Borst et al. \cite{crisfield2012nonlinear} (which however contains misprints regarding
the expressions for $\frac{\partial I_{2}}{\partial\mathbf{C}}$ and $\frac{\partial^{2}I_{3}}{\partial\mathbf{C}^{2}}$).

In Voigt notation, the deformation gradient is computed by

\[
\mathbf{F}=\begin{pmatrix}1 & 0 & 0 & 0 & 1 & 0 & 0 & 0 & 1\end{pmatrix}^{\rm T}+\mathrm{grad}\mathbf{U},
\]
where $\mathrm{grad}{\bf U}=\mathbf{B}\mathbf{a}$ and 

\[
\mathbf{B}=\begin{bmatrix}\frac{\partial\varphi_{1}}{\partial x_{1}} & 0 & 0 & \frac{\partial\varphi_{2}}{\partial x_{1}} & 0 & 0 & \cdots\\
\frac{\partial\varphi_{1}}{\partial x_{2}} & 0 & 0 & \frac{\partial\varphi_{2}}{\partial x_{2}} & 0 & 0 & \cdots\\
\frac{\partial\varphi_{1}}{\partial x_{3}} & 0 & 0 & \frac{\partial\varphi_{2}}{\partial x_{3}} & 0 & 0 & \cdots\\
0 & \frac{\partial\varphi_{1}}{\partial x_{1}} & 0 & 0 & \frac{\partial\varphi_{2}}{\partial x_{1}} & 0 & \cdots\\
0 & \frac{\partial\varphi_{1}}{\partial x_{2}} & 0 & 0 & \frac{\partial\varphi_{2}}{\partial x_{2}} & 0 & \cdots\\
0 & \frac{\partial\varphi_{1}}{\partial x_{3}} & 0 & 0 & \frac{\partial\varphi_{2}}{\partial x_{3}} & 0 & \cdots\\
0 & 0 & \frac{\partial\varphi_{1}}{\partial x_{1}} & 0 & 0 & \frac{\partial\varphi_{2}}{\partial x_{1}} & \cdots\\
0 & 0 & \frac{\partial\varphi_{1}}{\partial x_{2}} & 0 & 0 & \frac{\partial\varphi_{2}}{\partial x_{2}} & \cdots\\
0 & 0 & \frac{\partial\varphi_{1}}{\partial x_{3}} & 0 & 0 & \frac{\partial\varphi_{2}}{\partial x_{3}} & \cdots
\end{bmatrix}.
\]

The right Cauchy-Green tensor is then given by $\mathbf{C}=\mathbf{F}^{\rm T}\mathbf{F}$.
The derivatives of the invariants with respect to $\mathbf{C}$ are
given by

\[
\frac{\partial I_{1}}{\partial\mathbf{C}}=\begin{pmatrix}1\\
1\\
1\\
0\\
0\\
0
\end{pmatrix},\ \frac{\partial I_{2}}{\partial\mathbf{C}}=\begin{pmatrix}\text{{$C_{33}+C_{22}$}}\\
\text{{$C_{11}+C_{33}$}}\\
\text{{$C_{22}+C_{11}$}}\\
-C_{12}\\
-C_{23}\\
-C_{13}
\end{pmatrix}\ \frac{\partial I_{3}}{\partial\mathbf{C}}=\begin{pmatrix}C_{22}C_{33}-C_{23}^{2}\\
C_{33}C_{11}-C_{31}^{2}\\
C_{11}C_{22}-C_{12}^{2}\\
C_{23}C_{31}-C_{33}C_{12}\\
C_{31}C_{12}-C_{11}C_{23}\\
C_{12}C_{23}-C_{22}C_{31}
\end{pmatrix}.
\]

The derivatives of the modified invariants are given by

\[
\begin{cases}
\frac{\partial J_{1}}{\partial\mathbf{C}} & =I_{3}^{-1/3}\frac{\partial I_{1}}{\partial\mathbf{C}}-\frac{1}{3}I_{1}I_{3}^{-4/3}\frac{\partial I_{3}}{\partial\mathbf{C}}\\
\frac{\partial J_{2}}{\partial\mathbf{C}} & =I_{3}^{-2/3}\frac{\partial I_{2}}{\partial\mathbf{C}}-\frac{2}{3}I_{2}I_{3}^{-5/3}\frac{\partial I_{3}}{\partial\mathbf{C}}\\
\frac{\partial J_{3}}{\partial\mathbf{C}} & =\frac{1}{2}I_{3}^{-1/2}\frac{\partial I_{3}}{\partial\mathbf{C}}
\end{cases}
\]

The second derivatives of the invariants are given by

\[
\frac{\partial^{2}I_{1}}{\partial\mathbf{C}^{2}}=\mathbf{0},
\]

\[
\frac{\partial^{2}I_{2}}{\partial\mathbf{C}^{2}}=\begin{bmatrix}0 & 1 & 1 & 0 & 0 & 0\\
1 & 0 & 1 & 0 & 0 & 0\\
1 & 1 & 0 & 0 & 0 & 0\\
0 & 0 & 0 & -1/2 & 0 & 0\\
0 & 0 & 0 & 0 & -1/2 & 0\\
0 & 0 & 0 & 0 & 0 & -1/2
\end{bmatrix},
\]

and

\[
\frac{\partial^{2}I_{3}}{\partial\mathbf{C}^{2}}=\begin{bmatrix}0 & C_{33} & C_{22} & 0 & -C_{23} & 0\\
C_{33} & 0 & C_{11} & 0 & 0 & -C_{31}\\
C_{22} & C_{11} & 0 & -C_{12} & 0 & 0\\
0 & 0 & -C_{12} & -C_{33}/2 & C_{31}/2 & C_{23}/2\\
-C_{23} & 0 & 0 & C_{31}/2 & -C_{11}/2 & C_{12}/2\\
0 & \text{{$-C_{31}$} }& 0 & \text{{$-C_{23}/2$}} & C_{12}/2 & -C_{22}/2
\end{bmatrix} .
\]

The second derivatives of the modified invariants are then given by

\[
\begin{cases}
\frac{\partial^{2}J_{1}}{\partial\mathbf{C}^{2}} & =I_{3}^{-1/3}\frac{\partial^{2}I_{1}}{\partial\mathbf{C}^{2}}+\frac{4}{9}I_{1}I_{3}^{-7/3}\frac{\partial I_{3}}{\partial\mathbf{C}}\left(\frac{\partial I_{3}}{\partial\mathbf{C}}\right)^{\rm T}\\
 & \quad -\frac{1}{3}I_{3}^{-4/3}\left[\frac{\partial I_{1}}{\partial\mathbf{C}}\left(\frac{\partial I_{1}}{\partial\mathbf{C}}\right)^{\rm T}+I_{1}\frac{\partial^{2}I_{3}}{\partial\mathbf{C}^{2}}+\frac{\partial I_{3}}{\partial\mathbf{C}}\left(\frac{\partial I_{3}}{\partial\mathbf{C}}\right)^{\rm T}\right]\\
\frac{\partial^{2}J_{2}}{\partial\mathbf{C}^{2}} & =I_{3}^{-2/3}\frac{\partial^{2}I_{2}}{\partial\mathbf{C}^{2}}+\frac{10}{9}I_{2}I_{3}^{-8/3}\frac{\partial I_{3}}{\partial\mathbf{C}}\left(\frac{\partial I_{3}}{\partial\mathbf{C}}\right)^{\rm T}\\
 & \quad  -\frac{2}{3}I_{3}^{-5/3}\left[\frac{\partial I_{2}}{\partial\mathbf{C}}\left(\frac{\partial I_{3}}{\partial\mathbf{C}}\right)^{\rm T}+I_{2}\frac{\partial^{2}I_{3}}{\partial\mathbf{C}^{2}}+\frac{\partial I_{3}}{\partial\mathbf{C}}\left(\frac{\partial I_{2}}{\partial\mathbf{C}}\right)^{\rm T}\right]\\
\frac{\partial^{2}J_{3}}{\partial\mathbf{C}^{2}} & =\frac{1}{2}I_{3}^{-1/2}\frac{\partial^{2}I_{3}}{\partial\mathbf{C}^{2}}-\frac{1}{4}I_{3}^{-3/2}\frac{\partial I_{3}}{\partial\mathbf{C}}\left(\frac{\partial I_{3}}{\partial\mathbf{C}}\right)^{\rm T}
\end{cases}.
\]

The second Piola-Kirchhoff tensor on Voigt form is given as:

\[
\begin{cases}
\mathbf{S}_{\text{iso}} & =2\left(K_{1}\frac{\partial J_{1}}{\partial\mathbf{C}}+K_{2}\frac{\partial J_{2}}{\partial\mathbf{C}}\right)\\
\mathbf{S}_{\text{vol}} & =2K\left(J_{3}-1\right)\frac{\partial J_{3}}{\partial\mathbf{C}}
\end{cases}
\]
where $J_{3}=I_{3}^{\frac{1}{2}}$ . The constitutive forth order tensor is given on Voigt form as:

\[
\begin{cases}
\mathbf{L}_{\text{iso}} & =4\left(K_{1}\frac{\partial^{2}J_{1}}{\partial\mathbf{C}^{2}}+K_{2}\frac{\partial^{2}J_{2}}{\partial\mathbf{C}^{2}}\right)\\
\mathbf{L}_{\text{vol}} & =4\left(K\left(J_{3}-1\right)\frac{\partial^{2}J_{3}}{\partial\mathbf{C}^{2}}+K\frac{\partial J_{3}}{\partial\mathbf{C}}\otimes\frac{\partial J_{3}}{\partial\mathbf{C}}\right)
\end{cases}.
\]

The B-matrices used in the previous sections are given as follows:

\[
\mathbf{B}_{L}=\begin{bmatrix}F_{11}\frac{\partial\varphi^{1}}{\partial x} & F_{21}\frac{\partial\varphi^{1}}{\partial x} & F_{31}\frac{\partial\varphi^{1}}{\partial x} & \cdots & \cdots\\
F_{12}\frac{\partial\varphi^{1}}{\partial y} & F_{22}\frac{\partial\varphi^{1}}{\partial y} & F_{32}\frac{\partial\varphi^{1}}{\partial y} & \cdots & \cdots\\
F_{13}\frac{\partial\varphi^{1}}{\partial z} & F_{23}\frac{\partial\varphi^{1}}{\partial z} & F_{33}\frac{\partial\varphi^{1}}{\partial z} & \cdots & \cdots\\
F_{11}\frac{\partial\varphi^{1}}{\partial y}+F_{12}\frac{\partial\varphi^{1}}{\partial x} & F_{21}\frac{\partial\varphi^{1}}{\partial y}+F_{22}\frac{\partial\varphi^{1}}{\partial x} & F_{31}\frac{\partial\varphi^{1}}{\partial y}+F_{32}\frac{\partial\varphi^{1}}{\partial x} & \cdots & \cdots\\
F_{12}\frac{\partial\varphi^{1}}{\partial z}+F_{13}\frac{\partial\varphi^{1}}{\partial y} & F_{22}\frac{\partial\varphi^{1}}{\partial z}+F_{23}\frac{\partial\varphi^{1}}{\partial y} & F_{32}\frac{\partial\varphi^{1}}{\partial z}+F_{33}\frac{\partial\varphi^{1}}{\partial y} & \cdots & \cdots\\
F_{13}\frac{\partial\varphi^{1}}{\partial x}+F_{11}\frac{\partial\varphi^{1}}{\partial z} & F_{23}\frac{\partial\varphi^{1}}{\partial x}+F_{21}\frac{\partial\varphi^{1}}{\partial z} & F_{33}\frac{\partial\varphi^{1}}{\partial x}+F_{31}\frac{\partial\varphi^{1}}{\partial z} & \cdots & \cdots
\end{bmatrix},
\]

\[
\mathbf{B}_{NL}=\begin{bmatrix}\frac{\partial\varphi^{1}}{\partial x} & 0 & 0 & \cdots & \cdots\\
\frac{\partial\varphi^{1}}{\partial y} & 0 & 0 & \cdots & \cdots\\
\frac{\partial\varphi^{1}}{\partial z} & 0 & 0 & \cdots & \cdots\\
0 & \frac{\partial\varphi^{1}}{\partial x} & 0 & \cdots & \cdots\\
0 & \frac{\partial\varphi^{1}}{\partial y} & 0 & \cdots & \cdots\\
0 & \frac{\partial\varphi^{1}}{\partial z} & 0 & \cdots & \cdots\\
0 & 0 & \frac{\partial\varphi^{1}}{\partial x} & \cdots & \cdots\\
0 & 0 & \frac{\partial\varphi^{1}}{\partial y} & \cdots & \cdots\\
0 & 0 & \frac{\partial\varphi^{1}}{\partial z} & \cdots & \cdots
\end{bmatrix},
\]
and for the higher order terms in the stabilization we define the
derivatives of the B-matrices as follows:


\[
\frac{\partial\mathbf{B}_{L}}{\partial\xi_i}=\begin{bmatrix}F_{11}\varphi_{x}^{1} & F_{21}\varphi_{x}^{1} & F_{31}\varphi_{x}^{1} & \cdots & \cdots\\
F_{12}\varphi_{y}^{1} & F_{22}\varphi_{y}^{1} & F_{32}\varphi_{y}^{1} & \cdots & \cdots\\
F_{13}\varphi_{z}^{1} & F_{23}\varphi_{z}^{1} & F_{33}\varphi_{z}^{1} & \cdots & \cdots\\
F_{11}\varphi_{y}^{1}+F_{12}\varphi_{x}^{1} & F_{21}\varphi_{y}^{1}+F_{22}\varphi_{x}^{1} & F_{31}\varphi_{y}^{1}+F_{32}\varphi_{x}^{1} & \cdots & \cdots\\
F_{12}\varphi_{z}^{1}+F_{13}\varphi_{y}^{1} & F_{22}\varphi_{z}^{1}+F_{23}\varphi_{y}^{1} & F_{32}\varphi_{z}^{1}+F_{33}\varphi_{y}^{1} & \cdots & \cdots\\
F_{13}\varphi_{x}^{1}+F_{11}\varphi_{z}^{1} & F_{23}\varphi_{x}^{1}+F_{21}\varphi_{z}^{1} & F_{33}\varphi_{x}^{1}+F_{31}\varphi_{z}^{1} & \cdots & \cdots
\end{bmatrix},
\]
where $\varphi_{x}=\frac{\partial}{\partial\xi_{i}}\left(\frac{\partial\varphi}{\partial x}\right)$.
The derivatives
of $\mathbf{B}_{NL}$ are given as:

\[
\frac{\partial\mathbf{B}_{NL}}{\partial\xi_i}=\begin{bmatrix}\varphi_{x}^{1} & 0 & 0 & \varphi_{x}^{2} & 0 & 0 & \cdots\\
\varphi_{y}^{1} & 0 & 0 & \varphi_{y}^{2} & 0 & 0 & \cdots\\
\varphi_{z}^{1} & 0 & 0 & \varphi_{z}^{2} & 0 & 0 & \cdots\\
0 & \varphi_{x}^{1} & 0 & 0 & \varphi_{x}^{2} & 0 & \cdots\\
0 & \varphi_{y}^{1} & 0 & 0 & \varphi_{y}^{2} & 0 & \cdots\\
0 & \varphi_{z}^{1} & 0 & 0 & \varphi_{z}^{2} & 0 & \cdots\\
0 & 0 & \varphi_{x}^{1} & 0 & 0 & \varphi_{x}^{2} & \cdots\\
0 & 0 & \varphi_{y}^{1} & 0 & 0 & \varphi_{y}^{2} & \cdots\\
0 & 0 & \varphi_{z}^{1} & 0 & 0 & \varphi_{z}^{2} & \cdots
\end{bmatrix},
\]
where again $\varphi_{x}=\frac{\partial}{\partial\xi_{i}}\left(\frac{\partial\varphi}{\partial x}\right),$
etc.  Note that on
the tri-linear hexahedral element $\frac{\partial}{\partial\xi_{i}}\left(\frac{\partial\varphi}{\partial x}\right)$
is constant for i$=1,2,3$.

\bibliographystyle{elsarticle-num}
\bibliography{largedef}
\end{document}